\documentclass[]{interact}

\usepackage{epstopdf}
\usepackage[caption=false]{subfig}

\usepackage[numbers,sort&compress]{natbib}
\bibpunct[, ]{[}{]}{,}{n}{,}{,}
\makeatletter
\def\NAT@def@citea{\def\@citea{\NAT@separator}}
\makeatother
\usepackage{anyfontsize}
\theoremstyle{plain}
\newtheorem{theorem}{Theorem}[section]
\newtheorem{lemma}[theorem]{Lemma}
\newtheorem{corollary}[theorem]{Corollary}

\theoremstyle{definition}
\newtheorem{definition}[theorem]{Definition}
\newtheorem{remark}[theorem]{Remark}
\newtheorem{example}[theorem]{Example}

\usepackage{xcolor}

\theoremstyle{remark}

\begin{document}
	

	\title{Characterizations of closed EP operators on Hilbert spaces}
  \author{\name {Arup Majumdar \thanks{Arup Majumdar (corresponding author). Email address: arupmajumdar93@gmail.com}and P. Sam Johnson \thanks{ P. Sam Johnson. Email address: sam@nitk.edu.in}} \affil{Department of Mathematical and Computational Sciences, \\
			National Institute of Technology Karnataka, Surathkal, Mangaluru 575025, India.}}
  
	\maketitle

	\begin{abstract}
		In this paper, we present intriguing findings that characterize both the closed (unbounded) and bounded EP operators on Hilbert spaces. Additionally, we demonstrate the result $\gamma(T) \leq r(T)$, where $T$ is a bounded EP operator, and $\gamma(T) \text{ and } r(T)$ represent the reduced minimum modulus and the spectral radius of $T$, respectively.
	\end{abstract}
	
	\begin{keywords}
		EP operator, closed operator, Moore-Penrose inverse.
	\end{keywords}
      \begin{amscode}47A05; 47B02.\end{amscode}
    \section{Introduction}
  The term EP or Equal Projection was first introduced by Schwerdtfeger in 1950 \cite{MR001} to describe a square matrix $T$ over a complex field $\mathbb{C}$ for which the null spaces of $T$ and $T^{*}$ are identical. In 1965, M. H. Pearl established an equivalent condition of the EP matrix $T$ commutes with its Moore-Penrose inverse $T^{\dagger}$ \cite{MR002}. Subsequently, Campbell and Meyer expanded the concept of EP matrices to bounded linear operators on complex Hilbert spaces, defining EP operators as those for which the closed ranges of  $T$ and $T^{*}$ are equal \cite{MR003}. Itoh later introduced the hypo-EP operator, characterized by the conditions $T^{\dagger}T \geq TT^{\dagger}$ and $R(T)$ is closed \cite{MR005}. Research has continued into the characterization of EP operators on Hilbert spaces, with several authors investigating their properties within ${C^{*}}$-algebras. In 2007, Boasso studied EP operators in the context of  Banach space operators and Banach algebra elements \cite{MR004}, while Johnson focused on unbounded closed EP and hypo-EP operators on Hilbert spaces in 2021 \cite{MR006}.
  This paper delves into the exploration of the characterizations of closed (possibly unbounded) EP operators on Hilbert spaces, presenting intriguing examples in Section 2. Section 3 is dedicated to the examination of bounded EP operators on Hilbert spaces.

From now on, we shall mean $H$, $K$, $H_{i}$, $K_{i}$ ($i = 1,2,\dots, n$) as Hilbert spaces. The specification of a domain is an essential part of the definition of an unbounded operator, usually defined on a subspace.  Consequently, for an operator $T$, the specification of the subspace $D$ on which $T$ is defined, called the domain of $T$, denoted by $D(T)$, is to be given. The null space and range space of $T$ are denoted by $N(T)$ and $R(T)$, respectively. $W_{1}^{\perp}$ denotes the orthogonal complement of a set $W_{1}$ whereas $W_{1} \oplus^{\perp} W_{2}$ denotes the orthogonal direct sum of the subspaces $W_{1}$ and $W_{2}$ of a Hilbert space. Moreover,  $T_{\vert_{W}}$ denotes the restriction of $T$ to a subspace $W$ of a specified Hilbert space.  We call $D(T)\cap N(T)^\perp $, the carrier of $T$ and it is denoted by $ C(T)$. $T^{*}$ denotes the adjoint of $T$, when $D(T)$ is densely defined in the specified Hilbert space. Here, $P_{V}$ is the orthogonal projection on the closed subspace $V$ in the specified Hilbert space and the set of bounded operators from $H$ into $K$ is denoted by $B(H, K)$. Similarly, the set of all bounded operators on $H$ is denoted by $B(H)$.
	For the sake of completeness of exposition, we first begin with the definition of a closed operator.

	\begin{definition}
		Let $ T$ be an operator from a Hilbert space $H$ with domain $D(T) $ to a Hilbert space $K$. If the graph of $T$ defined by 
		$$ G(T)=\left\{(x,Tx): x\in D(T)\right\} $$ is closed in $H\times K $, then $T$ is called a closed  operator. Equivalently, $T $ is a closed operator if $ \{x_n \}$ in $D(T) $ such that $ x_n\rightarrow x $ and $ Tx_n\rightarrow y$ for some $ x\in  H,y\in   K $,  then $ x\in  D(T) $ and $ Tx=y $. That is, $G(T)$ is  a closed subspace of $H\times K$ with respect to the graph norm $\|(x, y)\|_T=(\|x\|^2+\|y\|^2)^{1/2}$. It is easy to show that the graph norm  $\|(x, y)\|_T$ is equivalent to the norm $\|x\|+\|y\|$. 	 We note that,  for any densely defined closed operator $T$,  the closure of $C(T)$, that is,  $ \overline{C(T)}$ is $ N(T)^\perp$. $C(H)$ denotes the set of all closed operators from $H$ into $H$. Whereas, $C(H, K)$ denotes the set of all closed operators from $H$ into $K$.
		
     We say that $S$ is an extension of $ T $ (denoted by $ T\subset S$) if $ D(T)\subset D(S) $ and $ Sx=Tx $ for all $ x\in D(T)$. Moreover, $T \in C(H)$ commutes $A \in B(H)$  if $AT \subset TA$.

	\end{definition}

	\begin{definition}
		Let $T$ be a closed operator from  $D(T) \subset H$ to  $K$. The Moore-Penrose inverse of $T$ is the map $T^{\dagger}: R(T) \oplus^{\perp} R(T)^{\perp} \to H$ defined by
		\begin{equation}\label{equ 1}
			T^{\dagger} y = 
			\begin{cases}
				({T}\vert_{C(T)})^{-1}y  & \text{if} ~  y\in R(T)\\
				0    & \text{if}  ~ y\in R(T)^{\perp}.
			\end{cases}
		\end{equation}
    \end{definition}
\noindent  It can be shown that $T^{\dagger}$ is bounded if and only if $R(T)$ is closed, when $T$ is closed.
	\begin{definition}\cite{MR003}\label{def 1.3}
 An operator $T \in B(H)$ is called an EP operator if $R(T)$ is closed and $R(T) = R(T^{*})$.
  \end{definition}
 \begin{definition}\cite{MR005}\label{def 1.4}
     An operator $T \in B(H)$ is called a hypo-EP operator if $R(T)$ is closed with $R(T) \subset R(T^{*})$.
 \end{definition}
 \begin{theorem}\cite{MR005}\label{thm 1.5}
 Let $T \in B(H, K)$ be a closed range operator and let $\{T_{n}\}$ be a sequence of closed range operators in $B(H, K)$. Let $T_{n}^{\dagger}$ be the Moore-Penrose inverse of $T_{n}$ for every $n$. Suppose that $T_{n} \to T$ (with respect to the norm of $B(H, K)$). Then the following conditions are equivalent:
 \begin{enumerate}
 \item $T_{n}^{\dagger} \to T$;
 \item $T_{n}^{\dagger}T_{n} \to T^{\dagger}T$;
 \item $\sup\{\| T_{n}^{\dagger}\| : n \in \mathbb{N}\} < \infty$.
 \end{enumerate}
 \end{theorem}
 Definitions \ref{def 1.3} and \ref{def 1.4} can be extended to densely defined closed operators, as discussed in the paper \cite{MR006}. Below, we present these definitions.
 \begin{definition}\cite{MR006}\label{def 1.6}
 Let $T$ be a densely defined closed operator on $H$. $T$ is said to be an EP operator if $T$ has a closed range and $R(T) =R(T^{*})$.\\
 Moreover, the densely defined closed operator $T \in C(H)$ is called a hypo-EP operator if $R(T)$ is closed with $R(T) \subset R(T^{*})$.
 \end{definition}

 \begin{theorem}\label{thm 1.7}\cite{MR0396607}
  Let $T$ be a densely defined closed operator from $D(T) \subset H$ into $K$. Then the following statements hold:
  \begin{enumerate}
  \item $T^{\dagger}$ is a closed operator from $K$ into $H$;
  \item $D(T^{\dagger}) = R(T) \oplus^{\perp} N(T^{*})$; $N(T^{\dagger}) = N(T^{*})$;
  \item $R(T^{\dagger}) = C(T)$;
  \item $T^{\dagger}Tx = P_{\overline{R(T^{\dagger})}}x, \text{ for all } x\in D(T)$;
  \item $TT^{\dagger}y = P_{\overline{R(T)}}y, \text{ for all } y\in D(T^{\dagger})$;
  \item $(T^{\dagger})^{\dagger} = T$;
  \item $(T^{*})^{\dagger} = (T^{\dagger})^{*}$;
  \item $N((T^{*})^{\dagger})= N(T)$;
  \item $(T^{*}T)^{\dagger} = T^{\dagger}(T^{*})^{\dagger}$;
  \item $(TT^{*})^{\dagger} = (T^{*})^{\dagger} T^{\dagger}$.
  \end{enumerate}
  \end{theorem}
 
 \begin{theorem}\label{thm 1.8} \cite{MR007}
Let $T_{1}: D(T_{1}) \subset H_{1} \to K_{1}$ and $T_{2}: D(T_{2}) \subset  H_{2} \to K_{2}$ be two closed operators with closed ranges. Then $T = T_{1} \bigoplus T_{2} : D(T_{1}) \bigoplus D(T_{2}) \subset H_{1} \bigoplus H_{2} \to K_{1} \bigoplus K_{2}$ has the Moore-Penrose inverse. Moreover,
 \begin{align*}
 T^{\dagger} = (T_{1} \bigoplus T_{2})^{\dagger} = T_{1}^{\dagger} \bigoplus T_{2}^{\dagger}.
 \end{align*}
\end{theorem}
\begin{theorem}\label{thm 1.9}\cite{MR007}
 Let $T_{i}:  D(T_{i}) \subset H_{i} \to K_{i} ~(i= 1, 2)$ be two closed operators with closed ranges $R(T_{i}), ~i = 1,2$. Then $\gamma(T_{1} \bigoplus T_{2}) = \min \{\gamma(T_{1}) ,\gamma(T_{2})\} > 0$, where $\gamma(T)$ is the reduced minimum modulus of $T$.
 \end{theorem}
 \section{Characterizations of closed EP operators on the Hilbert space $H$:}
 Throughout the section, we consider $T$ as the densely defined closed operator on $H$. Theorem \ref{thm 2.1} has been mentioned in the paper\cite{MR006} but we present a different approach to prove the theorem.
 \begin{theorem}\label{thm 2.1}
 Let $T \in C(H)$ be a closed-range operator. Then $T$ is EP if and only if $T^{\dagger}T = TT^{\dagger}$ on $D(T)$.
 \end{theorem}
 \begin{proof}
 Since, $R(T) = R(T^{*})$. Then, 
 \begin{align*}
 T^{\dagger}T = {P_{N(T)^{\perp}}}_{\vert{D(T)}} = {P_{R(T)}}_{\vert{D(T)}} \subset P_{R(T)} = TT^{\dagger}.
 \end{align*}
 Conversely, Suppose $T^{\dagger}T = TT^{\dagger}$ on $D(T)$. So,
 \begin{equation}\label{equ 2}
 R(T) = R(TT^{\dagger}) \supset \overline{R(T^{\dagger}T)} = \overline{R(T^{\dagger})}= N(T)^{\perp} = R(T^{*}).
  \end{equation}
  Again, $N(T) \subset R(T)^{\perp}$ because of $P_{R(T)}x = 0$, for all $x \in N(T)$. Thus,
  \begin{align}\label{equ 3}
  R(T) \subset N(T)^{\perp} = R(T^{*}).
  \end{align}
  By (\ref{equ 2}) and (\ref{equ 3}), we get $R(T) = R(T^{*})$.
 \end{proof}

 Now, we present some examples of EP operators.
\begin{example}\cite{MR006}
Let $\phi : [0,1] \to \mathbb{C}$ by
\begin{align*}
\phi(t) = 
\begin{cases}
1 \text{ if } t=0\\
\frac{1}{\sqrt{t}} \text{ if } 0< t \leq 1.
\end{cases}
\end{align*}
Define an operator $Tf = \phi f$ on $D(T) = \{f \in L^{2}[0,1] : \phi f \in L^{2}[0,1]\}$. $T$ is a self-adjoint operator. Moreover. $R(T)= L^{2}[0,1]$ because $\vert{\phi(t)\vert \geq 1}$. Therefore, $T$ is an EP operator.
\end{example}
\begin{example}\cite{majumdar2024hyers}
Define $T$ on $\ell^{2}$ by $$T(x_{1}, x_{2}, x_{3},\ldots,x_{n},\ldots) =(x_{1}, 2x_{2}, 3x_{3},\ldots,nx_{n},\ldots )$$ with domain $D(T)=\{(x_{1}, x_{2}, x_{3},\ldots,x_{n},\ldots)\in \ell^{2} : \sum_{n=1}^\infty |n x_{n}|^2<\infty\}$. Since  
		$D(T)$ contains the space $c_{00}$ of all finitely non-zero sequences, $ D(T)$ is a proper dense subspace of $\ell^{2}$. 
One can show that $T$ is a self-adjoint operator and  ${R(T)} = N(T^{*})^{\perp}= \ell^{2}= R(T^{*})$. Therefore, $T$ is EP.
\end{example}
\begin{example}\cite{majumdar2024hyers}
Define $T$ on  $\ell^{2}$ by $$T(x_{1}, x_{2},\ldots,x_{n},\ldots)= \Big(x_{1}, 2x_{2}, \frac{x_{3}}{3}, 4x_{4}, \frac{x_{5}}{5},\ldots\Big)$$ with domain $D(T)= \{(x_{1}, x_{2}, x_{3},\ldots,x_{n},\ldots)  : (x_{1}, 2x_{2}, \frac{x_{3}}{3}, 4x_{4}, \frac{x_{5}}{5},\ldots)\in \ell^{2} \}$. One can show that $T$ is closed and $R(T)$ is a proper dense subspace of $\ell^{2}$.  Therefore, $R(T)$ is not closed, which confirms that $T$ is not EP.
\end{example}
 \begin{theorem}\label{thm 2.2}
 Let $T_{1}\in C(H_{1})$ and $T_{2} \in C(H_{2})$ be two densely defined operators. Then $T_{1} \text{ and } T_{2}$ both are EP if and only if $T_{1} \bigoplus T_{2}$ is also EP.
 \end{theorem}
 \begin{proof}
 $R(T_{1})$ and $R(T_{2})$ both are closed because of EP operators. So, $R(T_{1} \bigoplus T_{2})$ is closed. From Theorem \ref{thm 1.8}, we have
 \begin{align*}
 (T_{1} \bigoplus T_{2})^{\dagger} (T_{1} \bigoplus T_{2}) &= T_{1}^{\dagger}T_{1} \bigoplus T_{2}^{\dagger}T_{2} \\ &= T_{1}T_{1}^{\dagger} \bigoplus T_{2}T_{2}^{\dagger} ~~~(\text{ on } D(T_{1}) \bigoplus D(T_{2}))\\ &= (T_{1} \bigoplus T_{2}) (T_{1} \bigoplus T_{2})^{\dagger} ~~~(\text{ on } D(T_{1} \bigoplus T_{2})).
 \end{align*}
 Hence, $(T_{1} \bigoplus T_{2})$ is an EP operator.\\
 Conversely, we claim that $R(T_{1})$ is closed. Consider an element $u \in \overline{R(T_{1})}$, there exists a sequence $\{T_{1}u_{n}\}$ in $R(T_{1})$ such that $T_{1}u_{n} \to u$, as $n \to \infty$. Since, $R(T_{1} \bigoplus T_{2})$ is closed. So, there is an element $(v_{1}, v_{2})$ in $D(T_{1} \bigoplus T_{2})$ such that
 $(T_{1} \bigoplus T_{2})(u_{n}, 0) \to  (u,0) = (T_{1}v_{1}, T_{2}v_{2})$, as $n \to \infty$ which implies $u = T_{1}v_{1}$ in $R(T_{1})$. Thus, $R(T_{1})$ is closed. Similarly, $R(T_{2})$ is closed.
 Now, for all $h_{1} \in D(T_{1})$, we get,
 \begin{align*}
 (T_{1}^{\dagger}T_{1}h_{1}, 0) &= (T_{1} \bigoplus T_{2})^{\dagger} (T_{1} \bigoplus T_{2})(h_{1}, 0)\\
 &= (T_{1} \bigoplus T_{2}) (T_{1} \bigoplus T_{2})^{\dagger} (h_{1},0)\\
 &=(T_{1}T_{1}^{\dagger}h_{1},0).
 \end{align*}
 This shows that $T_{1}^{\dagger}T_{1} = T_{1}T_{1}^{\dagger}$ on $D(T_{1})$. Similarly, it can be shown that $T_{2}^{\dagger}T_{2}= T_{2}T_{2}^{\dagger}$ on $D(T_{2})$. Therefore, $T_{1}$ and $T_{2}$ both are EP operators.
 \end{proof}
 \begin{theorem}\label{thm 2.3}
 Let $T \in C(H)$ be a closed range operator and $T= U_{T}\vert{T}\vert$ be the polar decomposition of $T$. Then the following conditions are equivalent:
 \begin{enumerate}
 \item $T$ is EP;
 \item $U_{T}$ is EP.
 \end{enumerate}
 \end{theorem}
 \begin{proof}
 $1) \Rightarrow 2)$ Since, $R(T) = R(T^{*})$. Then,
 \begin{align*}
 R(\vert T\vert) = R(T^{*}T) = R(T^{*}) = R(T).
 \end{align*}
 So, $R(U_{T}) = R(T)$ is closed and $N(U_{T}) = R(T)^{\perp} = N(T)$. Now, $U_{T}U_{T}^{\dagger} = P_{R(U_{T})} = P_{R(T)}$.
 Moreover, $U_{T}^{\dagger}U_{T} = P_{N(U_{T})^{\perp}} = P_{N(T)^{\perp}} = P_{R(T)}$. Thus, $U_{T}$ is EP.\\
 
 $2) \Rightarrow 1)$ Since, $U_{T}$ is EP. Then,
 \begin{align*}
 U_{T}U_{T}^{\dagger} &= U_{T}^{\dagger}U_{T}\\
 P_{R(U_{T})} &= P_{R(({U_{T}})^{*})}\\
 P_{R(T)} &= P_{R(T^{*})} ~~~~(\text{ because } R(U_{T}) = R(T) \text{ and } R((U_{T})^{*}) = R(U_{T^{*}}) = R(T^{*})).
 \end{align*}
 Hence, $R(T) = R(T^{*})$. Therefore, $T$ is EP.
 \end{proof}
 \begin{theorem}\label{thm 2.4}
 Let $T \in C(H)$ be of the form 
 \[
\left[\begin{array}{cc}
T_{1} & 0 \\
0 & 0
\end{array}\right]:\left[\begin{array}{c}
C(T) \\
N(T)
\end{array}\right] \rightarrow\left[\begin{array}{c}
R(T) \\
R(T)^{\perp}
\end{array}\right]
.\]\\
Then $T$ is EP if and only if $T_{1}$ is EP.
\end{theorem}
\begin{proof}
We claim that $T_{1} = T_{\vert{C(T)}}$ is EP. It is easy to show that $T_{1}$ is closed. Since $T$ is EP. So, $R(T) =R(T^{*})$. Again, $R(T_{1}) =R(T)$ is closed. Now, $T_{1}T_{1}^{\dagger} = I_{R(T)} = I_{N(T)^{\perp}}$ and $T_{1}^{\dagger}T_{1} = I_{C(T)}$. Thus, $T_{1}^{\dagger}T_{1} = T_{1}T_{1}^{\dagger}$ on $D(T_{1}) = C(T)$. Hence, $T_{1}$ is EP.\\
Conversely, $R(T) = R(T_{1})$ is closed. By Lemma 3.3 \cite{MR007}, It is obvious that $R(T^{*}) = R(T_{1}^{*})$. So, $T_{1}$ confirms that $R(T^{*}) = R(T_{1}^{*}) = R(T_{1}) = R(T)$. Therefore, $T$ is EP.
\end{proof}

\begin{theorem}\label{thm 2.5}
Let $T\in C(H)$ be an EP operator and $S$ be a bounded operator on $H$. Then, $T$ commutes with $S$ ($ST \subset TS$) if and only if $ST^{\dagger} = T^{\dagger}S$.
\end{theorem}
\begin{proof}
$R(T) = R(T^{*})$ is closed because $T$ is EP. So,
\[T=
\left[\begin{array}{cc}
 T_{1} & 0 \\
0 & 0
\end{array}\right]:\left[\begin{array}{c}
C(T) \\
N(T)
\end{array}\right] \rightarrow\left[\begin{array}{c}
R(T) \\
N(T)
\end{array}\right]
\] and 
\[T^{\dagger}=
\left[\begin{array}{cc}
 T_{1}^{-1} & 0 \\
0 & 0
\end{array}\right] \text{ is bounded}.\]
Moreover,
\[S=
\left[\begin{array}{cc}
S_{1} & S_{2} \\
S_{3} & S_{4}
\end{array}\right]:\left[\begin{array}{c}
R(T) \\
N(T)
\end{array}\right] \rightarrow\left[\begin{array}{c}
R(T) \\
N(T)
\end{array}\right]
.\]\\
For all $u =(u_{1}, u_{2}) \in D(ST)$, we have $STu = TSu$. Then,
\begin{align}\label{equ 4}
S_{1}T_{1}u_{1} = T_{1}(S_{1}u_{1} + S_{2}u_{2})
\end{align}
\begin{align}\label{equ 5}
S_{3}T_{1}u_{1} = 0.
\end{align}
Consider $w_{1} \in R(T)$, there exists an element $v_{1} \in C(T)$ such that $T_{1}^{-1}w_{1} = v_{1}$. From (\ref{equ 5}), we get $S_{3}w_{1} = S_{3}T_{1}v_{1}= 0$. Thus, $S_{3} = 0$. Choosing $u_{1} =0$ in (\ref{equ 4}), we get $T_{1}S_{2}u_{2} = 0$ which implies $S_{2} = 0$. Now, we claim that $S_{1}T_{1}^{-1} = T_{1}^{-1}S_{1}$.
When $u_{2} =0$, then for all $u_{1} \in C(T)$, (\ref{equ 4}) says that $S_{1}T_{1}u_{1} = T_{1}S_{1}u_{1}$.\\
For arbitrary $p_{1} \in R(T_{1})$, there is $z_{1} \in C(T)$ such that $z_{1} = T_{1}^{-1}p_{1}$. Again, $S_{1}p_{1} = T_{1}S_{1}T_{1}^{-1}p_{1}$. Thus, $T_{1}^{-1}S_{1}p_{1} = S_{1}T_{1}^{-1}p_{1}$. Hence, $T_{1}^{-1}S_{1} = S_{1}T_{1}^{-1}$. It is easy to show that
\[ST^{\dagger}=
\left[\begin{array}{cc}
S_{1}T_{1}^{-1} & 0 \\
0 & 0
\end{array}\right]\] and 
\[T^{\dagger}S=
\left[\begin{array}{cc}
T_{1}^{-1}S_{1} & 0 \\
0 & 0
\end{array}\right]\]
Therefore, $ST^{\dagger} = T^{\dagger}S$.\\
Conversely,  It is obvious that $T_{1}T_{1}^{-1}= I_{R(T)}$ and $T_{1}^{-1}T_{1}= I_{C(T)}$. Since, $ST^{\dagger} = T^{\dagger}S$. Then,
for all $x_{1} \in R(T) \text{ and } x_{2} \in N(T)$, 
\begin{align}\label{equ 6}
S_{1}T_{1}^{-1}x_{1} = T_{1}^{-1}(S_{1}x_{1} + S_{2}x_{2})
\end{align}
\begin{align}\label{equ 7}
S_{3}T_{1}^{-1}x_{1}=0
\end{align}
The equality (\ref{equ 7}) says that $S_{3} = 0$. Now consider, $x_{1}=0$, we get $S_{2}x_{2}=0$, for all $x_{2} \in N(T)$. Thus, $S_{2}=0$. Again taking, $x_{2}=0$, we have $S_{1}T_{1}^{-1}x_{1}= T_{1}^{-1}S_{1}x_{1}$ for all $x_{1}\in R(T)$ which implies $T_{1}S_{1}y = S_{1}T_{1}y$, for all $y\in C(T)$. So, for all $w_{1} \in C(T), w_{2} \in N(T)$,
\[ST\begin{pmatrix}
w_{1}\\
w_{2}
\end{pmatrix}=
\begin{pmatrix}
S_{1}T_{1}w_{1}\\
0
\end{pmatrix}=
\begin{pmatrix}
T_{1}S_{1}w_{1}\\
0
\end{pmatrix}=
TS\begin{pmatrix}
w_{1}\\
w_{2}
\end{pmatrix}.\]
Therefore, $ST \subset TS$.
\end{proof}
\begin{theorem}\label{thm 2.6}
Let $T \in C(H)$ be a closed-range operator. Then $T$ is EP if and only if $T^{n}$ is EP, for all $n\in \mathbb{N}$.~ (Here, we assume that $T^{n}$ is densely defined, for all $n \in \mathbb{N}$).
\end{theorem}
\begin{proof}
Suppose $T$ is EP. Then $R(T) = R(T^{*})$ is closed. So,
\begin{align}\label{equ 8}
R((T^{*})^{2}) = T^{*}R(T^{*})= T^{*}R(T)= R(T^{*})= R(T).
\end{align}
Again,
\begin{align}\label{equ 9}
R(T^{2}) = TR(T^{*})= R(TT^{*}) = R(T).
\end{align}
\begin{align}\label{equ 10}
R((T^{*})^{2}) \subset R((T^{2})^{*}) ~(\text{ because } (T^{*})^{2} \subset (T^{2})^{*}).
\end{align}
We know that $N(T) \subset N(T^{2})$. Thus,
\begin{align}\label{equ 11}
R((T^{2})^{*}) \subset N(T^{2})^{\perp} \subset N(T)^{\perp} = R(T^{*}) = R(T) = R((T^{*})^{2}).
\end{align}
By (\ref{equ 8}), (\ref{equ 9}), (\ref{equ 10}) and (\ref{equ 11}), we have $R((T^{2})^{*}) = R((T^{*})^{2}) = R(T) = R(T^{2})$. Since $R(T)$ is closed. Now we will show that $T^{2}$ is closed. Consider an element $(x, y) \in G(T^{2})$, where $G(T^{2})$ is the graph of $T^{2}$. Then there exists a sequence $\{x_{n}\}$ in $D(T^{2})$ such that $x_{n} \to x$ and $T^{2}x_{n} \to y$ as $n \to \infty$. The closed range of $T^{2}$ guarantees that there exists a element $z \in D(T^{2})$ with $y = T^{2}z$.  Again, the reduced minimal modulus of $T$, $\gamma(T) > 0$. So, from the relation $\|T^{2}x_{n} - T^{2}z\| \geq \gamma(T) \|Tx_{n}- Tz\|$, we get $Tx_{n} \to Tz$ as $n \to \infty$. Since $T$ is closed. So, $Tx = Tz$ which implies $x \in D(T^{2})$ and $y = T^{2}z = T^{2}x$. 

Hence, $T^{2}$ is EP. 
By induction hypothesis, it can be shown that $T^{n}$ is EP, for all $n \in \mathbb{N}$.
The converse part is obviously true.
\end{proof}
\begin{theorem}\label{thm 2.7}
Let $T \in C(H)$ be an EP operator. Then for all non zero $\lambda \in \mathbb{C}$, $\lambda \in \rho(T)$ if and only if $\lambda \in \rho(T_{\vert{C(T)}})$, where $T_{\vert{C(T)}} : C(T) \subset N(T)^{\perp} \to N(T)^{\perp} = R(T)$.
\end{theorem}
\begin{proof}
Suppose $\lambda \in \rho(T)$. Then $(T- \lambda)^{-1}$ is bounded. It is obvious that $N(T_{\vert{C(T)}} - \lambda) = \{0\}.$ Now consider $y \in N(T)^{\perp}$. There exists $x = x_{1} + x_{2} \in D(T)$, where $x_{1} \in C(T), x_{2}\in N(T)$, such that $y= (T- \lambda)x= Tx_{1} -\lambda x_{1}- \lambda x_{2}$. Thus, $y- (Tx_{1} -\lambda x_{1}) = -\lambda x_{2} \in N(T)^{\perp} \cap N(T) = \{0\}$. So, $y = (T_{\vert{C(T)}} - \lambda)x_{1}$. Thus, $T_{\vert{C(T)}}$ is onto. Again,
\begin{align*}
\|(T_{\vert{C(T)}} - \lambda)^{-1}y\| = \|x_{1}\| \leq \|x\| = \|(T-\lambda)^{-1}y\| \leq \|(T-\lambda)^{-1}\|\|y\|. 
\end{align*}
Hence, $(T_{\vert{C(T)}} - \lambda)^{-1}$ is bounded and $\lambda \in \rho(T_{\vert{C(T)}})$.\\
Conversely, Let $0 \neq \mu \in \rho(T_{\vert{C(T)}})$. Then $(T_{\vert{C(T)}} -\mu)^{-1}$ is bounded. We claim $(T- \mu)^{-1}$ is bounded on $H$. First consider, $z = z_{1} + z_{2} \in N(T- \mu)$, where $z_{1} \in C(T), z_{2} \in N(T)$. Then, $(T_{\vert{C(T)}} - \mu)z_{1} = \mu z_{2} \in N(T)^{\perp} \cap N(T) = \{0\}$. So, $z_{2} = 0$. The condition $N(T_{\vert{C(T)}} - \mu) = \{0\}$ says that $z_{1} = 0$. Hence, $N(T- \mu) = \{0\}$.\\
Now, we will show that $T - \mu$ is onto. Consider an arbitrary element $w = w_{1} + w_{2} \in H$, where $w_{1} \in R(T) = N(T)^{\perp} , w_{2} \in N(T)$. There is an element $u_{1} \in C(T)$ with $(T_{\vert{C(T)}} - \mu)u_{1} = w_{1}$. Taking the element $(u_{1} - \frac{w_{2}}{\mu}) \in D(T)$, we get 
\begin{align*}
(T- \mu)(u_{1} - \frac{w_{2}}{\mu}) = (T_{\vert{C(T)}} - \mu)u_{1} + w_{2} = w_{1} + w_{2} = w.
\end{align*}
This confirms that $T-\mu$ is onto. Moreover,
\begin{align*}
\|(T- \mu)^{-1}w\| = \|u_{1} - \frac{w_{2}}{\mu}\| \leq \|(T_{\vert{C(T)}} -\mu)^{-1}w_{1}\| + \frac{\|w_{2}\|}{|\mu|} &\leq M(\|w_{1}\| + \|w_{2}\|)\\
&\leq M \sqrt{2} \|w_{1} + w_{2}\|\\
&= M \sqrt{2} \|w\|
\end{align*}
, where $M = max\{ \|(T_{\vert{C(T)}} - \mu)^{-1}\|, \frac{1}{|\mu|}\}$.
Therefore, $\mu \in \rho(T)$.
\end{proof}
\begin{corollary}\label{cor 2.8}
Let $T \in C(H)$ be an EP operator. Then $0$ is not a limit point of the spectrum of $T$.
\end{corollary}
\begin{proof}
Since $T$ is EP. Then, $T_{\vert{C(T)}}$ is closed, one-one and onto, where $T_{\vert{C(T)}} : C(T) \to R(T)$. For all $y \in R(T)$, we have $\|(T_{\vert{C(T)}})^{-1}y\| = \|T^{\dagger}y\| \leq \|T^{\dagger}\| \|y\|$. So, $0 \in \rho(T_{\vert{C(T)}})$. We know that $\rho(T_{\vert{C(T)}})$ is open then there exists $\epsilon > 0$ such that $\{\lambda \in \rho(T_{\vert{C(T)}}): 0 < |\lambda| < \epsilon\} \subset \rho(T)$ ~(by Theorem \ref{thm 2.7}). Therefore, $0$ is not a limit point of the spectrum of $T$.  
\end{proof}
\begin{corollary}\label{cor 2.9}
Let $T \in C(H)$ be an EP operator with $\dim{N(T)} < \infty$. Then $0\in \Delta_{k}(T)$, for all $k \in \{1, 2, 3, 4,5\}$, where $\Delta_{k}(T) = \mathbb{C} \setminus \sigma_{ek}(T)$ and $\sigma_{ek}(T)$ are the various essential spectra defining in \cite{MR008}.
\end{corollary}
\begin{proof}
Since, $R(T)$ is closed and $R(T)^{\perp} = N(T)$. Then $0 \in \Delta_{k}(T)$, for all $k= 1,2,3,4$. Moreover, $\Delta_{5}(T)$ is the union of all the components of $\Delta_{1}(T)$ which intersect the resolvent set $\rho(T)$ of $T$. Therefore, by Corollary \ref{cor 2.8}, we get $0 \in \Delta_{5}(T)$.
\end{proof}
We consider two operators $S \in B(H)$ and $T \in C(H)$ but $ST$ is not closed in general. We illustrate an example to show that $ST$ is not closed.
\begin{example}\label{exam 2.10}
Let $S$ is defined on $\ell^{2}$ by:\\
\begin{align*}
S(x_{1}, x_{2}, \dots, x_{n},\dots) = (x_{1},{\frac{1}{2}}x_{2},\dots, {\frac{1}{n}}x_{n},\dots)
\end{align*}
Then, $S$ is bounded. Consider $T$ on $\ell^{2}$ by:
\begin{align*}
T(x_{1}, x_{2}, \dots, x_{n},\dots) = (x_{1}, 2x_{2},\dots, n x_{n},\dots)
\end{align*}
$T$ is self-adjoint. So, $T$ is closed. But $D(ST)= D(T) = \{x\in \ell^{2}: Tx \in \ell^{2}\}$ is densely defined in $\ell^{2}$. $ST= I_{D(T)}$ is bounded. If $ST$ is closed then $D(T) = \ell^{2}$, which is a contradiction. Therefore, $ST$ is not closed. 
\end{example}
The next Lemma \ref{lemma 2.11} gives a sufficient condition to have the closed operator $ST$ when $S \in B(H)$ and $T \in C(H)$. 
\begin{lemma}\label{lemma 2.11}
Let $T \in C(H)$ and $S \in B(H)$ be an EP operator with $R(ST) = R(T)$.  Then $ST$ is closed.
\end{lemma}
\begin{proof}
We can write $T$ and $S$ as following:
\[T=
\left[\begin{array}{cc}
 T_{1} & 0 \\
T_{2} & 0
\end{array}\right]:\left[\begin{array}{c}
C(T) \\
N(T)
\end{array}\right] \rightarrow\left[\begin{array}{c}
R(S) \\
N(S)
\end{array}\right]
\] and
\[S=
\left[\begin{array}{cc}
 S_{1} & 0 \\
0 & 0
\end{array}\right]:\left[\begin{array}{c}
R(S) \\
N(S)
\end{array}\right] \rightarrow\left[\begin{array}{c}
R(S) \\
N(S)
\end{array}\right]
.\]
It is obvious to show that $S_{1}^{-1}$ is bounded. From the given condition $R(ST) = R(T)$, we get $T_{2} = 0$. Thus, $T_{1}$ is closed. So, $S_{1}T_{1}$ is also closed. Now, consider an arbitrary element $(p,q) \in \overline{G(ST)}$. Then there is a sequence $\{p_{n}\} = \{p_{n}^{'} + p_{n}^{''}\}$ in $D(T)$ with $p_{n}^{'} \in C(T),~ p_{n}^{''} \in N(T)$, for all $n \in \mathbb{N}$, such that $p_{n}^{'} \to p^{'}$ and $p_{n}^{''} \to p^{''}$ as $ n \to \infty$, where $p = p^{'} + p^{''}, p^{'} \in R(T), \text{ and } p^{''}\in N(T)$. $ST{p_{n}} \to q = q^{'} +q^{''}$ as $n \to \infty$, where $q^{'} \in R(S), q^{''} \in N(S)$. So, $q^{''} =0$ and $S_{1}T_{1}p_{n}^{'} \to q^{'}$ as $n \to \infty$. The closedness of $S_{1}T_{1}$ confirms that $S_{1}T_{1}p^{'} = q^{'}$. Thus, 
\[q=
\begin{pmatrix}
S_{1}T_{1}p^{'}\\
0
\end{pmatrix}
= ST
\begin{pmatrix}
p^{'}\\
p^{''}
\end{pmatrix}
= STp.\]
Therefore, $ST$ is closed.
\end{proof}
If $A$ and $B$ are two EP operators, it remains an open question under what conditions the product $AB$ is also EP. In \cite{MR1428641}, Hartwig and Katz provided the necessary and sufficient conditions for the product of two square EP matrices to be EP. Subsequently, Dragan S. Djordjević established the necessary and sufficient conditions for the product of two bounded EP operators on $B(H)$ to also be EP \cite{MR1814103}. In Theorem \ref{thm 2.12}, we will present the necessary and sufficient conditions for the product of a closed EP operator and a bounded EP operator to also be EP.
\begin{theorem}\label{thm 2.12}
Let $T \in C(H)$ and $S\in B(H)$ both be EP with having $T^{*}$ has a matrix representation. Then, $ST$ is EP if and only if
\begin{enumerate}
\item R(ST) = R(T); 
\item N(ST) = N(T).
\end{enumerate}
\end{theorem}
\begin{proof}
Since $ST$ is EP. Then, $ST$ is closed and $R(ST) = R(T^{*}S^{*}) \subset R(T^{*}) = R(T)$.
\[\text{ For } \begin{pmatrix}
u\\
v
\end{pmatrix} \in \begin{bmatrix}
C(T)\\
N(T)
\end{bmatrix}, \text{ there exists $\begin{pmatrix}
w\\
t
\end{pmatrix} \in \begin{bmatrix}
R(S)\\
N(S)
\end{bmatrix}$ such that, }
\] 
\[
\left[\begin{array}{cc}
 T_{1} & 0 \\
T_{2} & 0
\end{array}\right]\begin{pmatrix}
u \\
v
\end{pmatrix}  =
\left[\begin{array}{cc}
 T_{1}^{*} & T_{2}^{*} \\
0 & 0
\end{array}\right]\begin{pmatrix}
w \\
t
\end{pmatrix}
.\] Thus,
\begin{align}
T_{1}u= T_{1}^{*}w + T_{2}^{*}t \text{ and }
T_{2}u= 0
\end{align}
So, $T_{2} = 0$ and $R(T) =R(T_{1}) \subset R(S)$. Moreover, $N(T_{1}) =\{0\}$. Consider $w_{1} + w_{2} \in N(ST)$, where $w_{1} \in C(T), w_{2} \in N(T)$, we get $S_{1}T_{1}w_{1} = 0$ which implies $w_{1} =0$. Thus, $ N(ST) \subset N(T) \subset N(ST)$. Thus, the condition (2), $N(ST) = N(T)$ is true. Taking an element $y \in R(T) \cap R(ST)^{\perp}$. Then, $y \in R(ST)^{\perp} = N(ST) = N(T)= N(T^{*}) = R(T)^{\perp}$ which implies $y= 0$. Hence, 
$R(T) \subset R(ST) \subset R(T)$. Therefore, the condition (1), $R(ST) =R(T)$ is also satisfied.\\
Conversely, Lemma \ref{lemma 2.11} confirms that $ST$ is closed. $R(ST) = R(T)$ is closed. Now,
\begin{align*}
R((ST)^{*}) &= N(ST)^{\perp}\\
&= N(T)^{\perp} ~(\text{ by condition (2)})\\
&= R(T)\\
&= R(ST) ~ ( \text{ by condition (1)}).
\end{align*}
It shows that  $ST$ is EP.
\end{proof}
\begin{theorem}\label{thm 2.13}
Let $T \in C(H, K)$ be a closed range operator. Then $R(\vert{T}\vert) = R({\vert{T}\vert}^{\alpha})$ for all $\alpha \in (0, \infty)$.
\end{theorem}
\begin{proof}
Since $R(T)$ is closed. Then $R(T^{*}) = R(T^{*}T) = R(\vert{T}\vert)$ is closed. Now,
\begin{align}\label{equ 13}
R(\vert{T}\vert) = R({\vert{T}\vert}^{\frac{1}{2}}) = \dots = R({\vert{T}\vert}^{\frac{1}{2^{n}}}), \text{ for all } n\in \mathbb{N}.
\end{align}
When $1 < m \leq 2^{n}$ and $m \in \mathbb{N}$, we have
\begin{align}
R({\vert{T}\vert}^{\frac{m}{2^{n}}}) &= {\vert{T}\vert}^{\frac{m-1}{2^{n}}} R({\vert{T}\vert}^{\frac{1}{2^{n}}})\\
&= {\vert{T}\vert}^{\frac{m-1}{2^{n}}} R({\vert{T}\vert}^{\frac{1}{2^{n-1}}})\\
&= R({\vert{T}\vert}^{\frac{m+1}{2^{n}}}).
\end{align}
By induction hypothesis, we can say that 
\begin{align}
R({\vert{T}\vert}^{\frac{m+k}{2^{n}}}) = R({\vert{T}\vert}^{\frac{m}{2^{n}}}), \text{ for all} ~k\in \mathbb{N}.
\end{align}
Let us choose $k = 2^{n} - m$, then $R({\vert{T}\vert}^{\frac{m}{2^{n}}}) = R(\vert{T}\vert)$, for all $m,n \in \mathbb{N}$ and $1 \leq m \leq 2^{n}$.\\
For all $\alpha \in (0,1)$, we have $\vert{T}\vert = {\vert{T}\vert}^{\alpha} {\vert{T}\vert}^{1-\alpha}$. So, 
\begin{align}\label{equ 18}
R(\vert{T}\vert) \subset R({\vert{T}\vert}^{\alpha}).
\end{align}
Again, we also have $m^{'}, n^{'} \in \mathbb{N}$ such that $1 \leq m^{'} \leq 2^{n^{'}}$ and $\frac{m^{'}}{2^{n^{'}}} \leq \alpha$.
Thus,
\begin{align}\label{equ 19}
R({\vert{T}\vert}^{\alpha}) \subset R({\vert{T}\vert}^{\frac{m^{'}}{2^{n^{'}}}}) = R({\vert{T}\vert}).
\end{align}
By (\ref{equ 18}) and (\ref{equ 19}), $R({\vert{T}\vert}^{\alpha}) = R({\vert{T}\vert})$, for all $\alpha \in (0,1]$.
We know that $R({\vert{T}\vert}^{2}) = R({\vert{T}\vert})$.
Let us assume that $R({\vert{T}\vert}^{q}) = R({\vert{T}\vert})$, where $q \geq 2$.
\begin{align*}
R({\vert{T}\vert}^{q+1}) = {\vert{T}\vert}^{q-1} R({\vert{T}\vert}^{2}) = {\vert{T}\vert}^{q-1} R({\vert{T}\vert}) =  R({\vert{T}\vert}^{q}) = R({\vert{T}\vert}).
\end{align*}
Then, $R({\vert{T}\vert}^{n}) = R({\vert{T}\vert})$, for all $n \in \mathbb{N}$.\\
Consider an arbitrary element $1 < \beta < \infty$, there exists $\mu \in \mathbb{N}$ such that $\mu < \beta \leq \mu + 1$. So, $0 < \beta - \mu \leq 1 $. Hence, 
$R({\vert{T}\vert}^{\beta}) = {\vert{T}\vert}^{\mu} R({\vert{T}\vert}^{\beta- \mu}) = R({\vert{T}\vert}^{\mu+1}) = R({\vert{T}\vert}).$\\
Therefore, $R(\vert{T}\vert) = R({\vert{T}\vert}^{\alpha})$, for all $\alpha \in (0, \infty).$
\end{proof}
\begin{corollary}\label{cor 2.14}
Let $T \in C(H)$ be an EP operator. Then, $R(T) = R(\vert{T}\vert)= R({\vert{T}\vert}^{\alpha})$, for all $\alpha \in (0, \infty)$.
\end{corollary}
\begin{proof}
Since $T$ is EP. By Theorem \ref{thm 2.13}, $ R({\vert{T}\vert}^{\alpha}) = R(\vert{T}\vert)= R(T^{*}T) = R(T^{*}) = R(T).$
\end{proof}
The following Theorem \ref{thm 2.15} says that the converse of Corollary \ref{cor 2.14} is true when $T \in B(H)$ and for some $\alpha \in (0,1) $.
\begin{theorem}\label{thm 2.15}
Let $T \in B(H)$ be satisfied the condition $R(T) = R(\vert{T}\vert) = R({\vert{T}\vert}^{\alpha})$, for some $\alpha \in (0,1)$. Then $T$ is EP.
\end{theorem}
\begin{proof}
First, we claim that $R(\vert{T}\vert)$ is closed. Since $T$ is bounded. So, $H = D({\vert{T}\vert}^{\alpha}) + R(\vert{T}\vert)$. For given $x \in N(\vert{T}\vert)^{\perp}$, there are $x_{1} \in D({\vert{T}\vert}^{\alpha}), x_{2} \in R(\vert{T}\vert)$ such that $x = x_{1} + x_{2}$. Then, $x_{1} \in N(\vert{T}\vert)^{\perp}$. For $x_{1}$ there exists $y$ such that $\vert{T}\vert^{\alpha}x_{1} = \vert{T}\vert y$ which implies $(x_{1} - \vert{T}\vert^{1-\alpha}y) \in N(\vert{T}\vert^{\alpha}) = N(\vert{T}\vert)$. Again, $\vert{T}\vert^{1-\alpha}y \in R(\vert{T}\vert^{1-\alpha}) \subset N(\vert{T}\vert^{1-\alpha})^{\perp} = N(\vert{T}\vert)^{\perp}$ (by McCarthy inequalities $N({\vert{T}\vert}^{\beta}) = N(\vert{T}\vert), \text{ for all $\beta \in (0,1]$}$). Thus, $(x_{1} - \vert{T}\vert^{1-\alpha}y) \in N(\vert{T}\vert) \cap N(\vert{T}\vert)^{\perp} = \{0\}$ which confirms that $x_{1} \in R({\vert{T}\vert}^{1-\alpha})$. Moreover, $x = x_{1} + x_{2} \in R({\vert{T}\vert}^{1-\alpha})$. Hence, 
\begin{align}
\overline{R({\vert{T}\vert}^{1-\alpha})} = \overline{R({\vert{T}\vert})} = N({\vert{T}\vert})^{\perp} \subset R({\vert{T}\vert}^{1-\alpha}) \subset \overline{R({\vert{T}\vert}^{1-\alpha})}.
\end{align}
So, $R({\vert{T}\vert}^{1-\alpha})$ is closed which implies $R({\vert{T}\vert}) \subset R({\vert{T}\vert}^{1-\alpha}) = R({\vert{T}\vert}^{2^{k}(1-\alpha)}) \subset R({\vert{T}\vert})$, for some $k \in \mathbb{N}$. Hence, $R({\vert{T}\vert}) = R(T)$ is closed. It means that $R(T^{*}) = R({\vert{T}\vert}^{2}) = R(\vert{T}\vert)= R(T)$. Therefore, $T$ is EP.
\end{proof}
\begin{theorem}\label{thm 2.16}
Let $T \in C(H)$ be an EP operator and $S \in B(H)$ be also an EP operator with $\|Sx\| \leq a\|Tx\|$, for all $x \in D(T)$ and $0 < a < 1$. Then $T + S$ is hypo-EP.
\end{theorem}
\begin{proof}
It is easy to show that $T + S$ is closed. The given condition says that
\begin{align}\label{equ 21}
(1-a) \|Tx\| \leq \|(T + S)x\| \leq (1 + a) \|Tx\|, \text{ for all}~ x\in D(T).
\end{align}
The above inequality (\ref{equ 21}) guarantees that $R(T + S)$ is closed. So, $R((T + S)^{*})$ is also closed. $N(T) \subset N(S)$ says that $N(T) \subset N(T + S)$. For all $u \in D(T^{*}T) \subset D(T)$, we get
\begin{align*}
\|Su\|^{2} &\leq a^{2} \|Tu\|^{2}\\
\langle S^{*}Su, u \rangle &\leq a^{2} \langle T^{*}Tu, u \rangle\\
S^{*}S &\leq a^{2}T^{*}T.
\end{align*}
By Douglas Theorem \cite{MR0203464}, there exists a contraction $C$ such that $S^{*} \subset (aT)^{*}C$. So, $R(S^{*}) \subset R(T^{*})$.
When an arbitrary $p \in N(T +S)$, then $\|Tp\| = \|Sp\| \leq a \|Tp\| < \|Tp\|$ confirms that $p \in N(T)$. Thus, $N(T + S) = N(T)$. Again, $R(T + S)^{*} = R(T^{*}) = R(T^{*}) + R(S^{*}) = R(T) + R(S) \supset R(T + S)$. Therefore, $T + S$ is hypo-EP.
\end{proof}
\begin{corollary}\label{cor 2.17}
Let $T\in C(H)$ and $S \in B(H)$ be EP with the following conditions:
\begin{enumerate}
\item $\|Sx\| \leq a \|Tx\|, \text{ where } a < 1 \text{ and for all } x \in D(T)$;
\item $\|S^{*}z\| \leq b \|T^{*}z\|, \text{ where } b < 1 \text{ and for all } z \in D(T^{*})$.
\end{enumerate}
Then $T + S$ is EP.
\end{corollary}
\begin{proof}
From Theorem \ref{thm 2.16}, we have $R((T +S)^{*}) = R(T)$ is closed and $T + S$ is closed. Again, the mentioned condition (2) guarantees that $R(T + S) = R(T)$. Thus, $R(T + S)= R((T + S)^{*}) = R(T)$. Therefore, $T + S$ is EP.
\end{proof}
\begin{remark}\label{remark 2.18}
Let $T \in C(H)$ and $S \in B(H)$ be two normal and EP operators with $\|Sx\| < a \|Tx\|, \text{ where } a < 1 \text{ and for all } x \in D(T)$. Then, $T +S$ is EP.
\end{remark}
\begin{theorem}\label{thm 2.19}
Let $T \in C(H)$ be a closed-range operator. Then $T$ is EP if and only if 
\begin{enumerate}
\item $T(I - TT^{\dagger})x = 0, \text{ for all } x \in H$;
\item $T^{*}(I - T^{*}(T^{*})^{\dagger})x = 0, \text{ for all } x \in H$.
\end{enumerate}
\end{theorem}
\begin{proof}
Suppose $T$ is EP. Then $N(T) = N(T^{*})$. For all $x \in H$, we have
\begin{align*}
T(I - TT^{\dagger})x = TP_{R(T)^{\perp}}x = TP_{N(T)}x = 0,
\end{align*}
and
\begin{align*}
T^{*}(I - T^{*}(T^{*})^{\dagger})x = T^{*}P_{N(T)}x = T^{*}P_{N(T^{*})}x=0.
\end{align*}
Conversely, From the condition $(1)$, we get $TP_{R(T)^{\perp}}x=0 , \text{ for all } x \in H$. So, $R(T)^{\perp} = N(T^{*}) \subset N(T)$. Condition$(2)$ says that $R(T^{*})^{\perp} = N(T) \subset N(T^{*})$. Therefore, $T$ is EP.
\end{proof}

\section{Characterizations of bounded EP operators on Hilbert spaces:} 
 
Let us consider the set $E$ of all EP operators in $B(H)$. The following example shows that $E$ is not closed in $B(H)$.
\begin{example}
Let us define $T_{n}: \ell^{2} \to \ell^{2}$ by
\begin{align*}
T_{n}(x_{1}, x_{2},\dots,x_{k},\dots) = (x_{1}, \frac{1}{2}x_{2},\dots, \frac{1}{n}x_{n},0,\dots,0,\dots), \text{ for all } n\in \mathbb{N}.
\end{align*}
Here, each $T_{n} \in B(\ell^{2})$ is self-adjoint closed range operator. So, $T_{n}$ is EP, for all $n \in \mathbb{N}$. Moreover $T_{n} \to T$, as $n \to \infty$, where $T: \ell^{2} \to \ell^{2}$ defined by:
\begin{align*}
T(x_{1}, x_{2},\dots,x_{k},x_{k+1}\dots) = (x_{1}, \frac{1}{2}x_{2},\dots, \frac{1}{k}x_{k},\frac{1}{k+1}x_{k+1},\dots). 
\end{align*}
It is easy to show that $\gamma(T) = 0$, where $\gamma(T)$ is the reduced minimum modulus of $T$.
Thus, $T$ is not EP because $R(T) = R(T^{*})$ but $R(T)$ is not closed.
\end{example}
The following Theorem \ref{thm 3.2} says that there is a closed set in $B(H)$ whose all elements are EP. Before stating Theorem \ref{thm 3.2}, a new set $E_{\delta}$ is  defined by:
\begin{align}
E_{\delta} = \{ T \in E  : \gamma(T) \geq \delta > 0\}.
\end{align}
\begin{theorem}\label{thm 3.2}
$E_{\delta}$ is closed set in $B(H)$.
\end{theorem}
\begin{proof}
Let us consider $T \in \overline{E_{\delta}}, \text{ the closure of $E_{\delta}$}$. Then, there is a sequence $\{T_{n}\}$ in $E_{\delta}$ such that $T_{n} \to T$ as $n \to \infty$. So, $\gamma(T_{n}) \geq \delta$ implies $\gamma(T) \geq \delta$ \cite{MR1888031}. Thus, $R(T)$ is closed. We know,$ \gamma(T_{n}) = \frac{1}{\|{T_{n}}^{\dagger}\|}, \text{ for all } n \in \mathbb{N}$. Thus, $\sup\{\|{T_{n}}^{\dagger}\|:n \in \mathbb{N}\} \leq \frac{1}{\delta}$. By Theorem \ref{thm 1.5}, we get ${T_{n}}^{\dagger} \to T^{\dagger}$ as $n \to \infty$. Hence,
\begin{align*}
&\|T^{\dagger}T - TT^{\dagger}\| \\
&= \|T^{\dagger}(T-T_{n})+ T^{\dagger}T_{n}-{T_{n}}^{\dagger}T_{n} + T_{n}{T_{n}}^{\dagger} -T_{n}T^{\dagger} + T_{n}T^{\dagger}-TT^{\dagger}\|\\
&\leq \|T^{\dagger}\|\|T-T_{n}\| + \|T^{\dagger}-{T_{n}}^{\dagger}\|\|T_{n}\|+ \|T_{n}\|\|{T_{n}}^{\dagger}-T^{\dagger}\|+ \|T_{n}-T\|\|T^{\dagger}\|.
\end{align*}
The right-hand side of the above inequality goes to $0$, as $n \to \infty$. Hence, $T$ is EP. Therefore, $E_{\delta}$ is closed in $B(H)$.
\end{proof}
\begin{corollary}
    $E = (\bigcup\limits_{\delta > 0} E_{\delta}) \bigcup \{0\}$.
\end{corollary}
\begin{proof}
Let $S \in E$ be a non-zero EP operator. Then, $\gamma(S) \geq \delta_{1}$, for some $\delta_{1} > 0$. So, $S \in E_{\delta_{1}}$. Thus, $E \subset (\bigcup\limits_{\delta > 0} E_{\delta}) \bigcup \{0\}$. The opposite inclusion is obvious. Therefore, $E = (\bigcup\limits_{\delta > 0} E_{\delta}) \bigcup \{0\}$.
\end{proof}

\begin{theorem}\label{thm 3.4}
Let $T \in B(H)$ be an EP operator. Then 
\begin{align*}
\gamma(T) \leq r(T).
\end{align*}
\end{theorem} 
\begin{proof}
Theorem \ref{thm 2.6} says that $T^{n}$ is EP and $R(T^{n}) = R(T)$, for all $n \in \mathbb{N}$. Moreover, $(T^{n})^{\dagger} = (T^{\dagger})^{n}$, for all $n\in \mathbb{N}$. We know that $r(T^{\dagger}) = \lim_{{n \to \infty}}\|(T^{\dagger})^{n}\|^{\frac{1}{n}}$. Again, 
\begin{align}\label{equ 23}
\frac{1}{\|(T^{\dagger})^{n}\|}=\frac{1}{\|(T^{n})^{\dagger}\|} = \gamma(T^{n}) = \displaystyle{\inf_{{x \in R(T^{n})}}} \frac{\|T^{n}x\|}{\|x\|} \leq (\displaystyle{\inf_{{x \in R(T)}}} \frac{\|Tx\|}{\|x\|}) \|T\|^{n-1} = \|T\|^{n-1} \gamma(T).
\end{align}
Thus, $\frac{1}{\|T\|}={\lim_{n \to \infty}}{\frac{1}{\|T\|^{1- \frac{1}{n}}}} \leq r(T^{\dagger})$. Therefore, $\gamma(T) = \frac{1}{\|T^{\dagger}\|} \leq r(T)$.
\end{proof}

\begin{center}
	\textbf{Acknowledgements}
\end{center}

\noindent The present work of the second author was partially supported by Science and Engineering Research Board (SERB), Department of Science and Technology, Government of India (Reference Number: MTR/2023/000471) under the scheme ``Mathematical Research Impact Centric Support (MATRICS)''. 

\section*{Declarations}
  The authors declare no conflicts of interest.

\end{document}